\documentclass{article}
\usepackage[a4paper, portrait, margin=1.1811in]{geometry}
\usepackage[english]{babel}
\usepackage[utf8]{inputenc}
\usepackage[T1]{fontenc}
\usepackage{helvet}
\usepackage{etoolbox}
\usepackage{graphicx}
\usepackage{titlesec}
\usepackage{caption}
\usepackage{booktabs}
\usepackage{xcolor} 
\usepackage{caption}
\graphicspath{ {./images/} }
\usepackage{scrextend}
\usepackage{fancyhdr}
\usepackage{graphicx}

\usepackage{hyperref}
\hypersetup{
   colorlinks=true,
   linkcolor=black,
    filecolor=black,      
    urlcolor=black,
    citecolor=black
    }

\usepackage{amssymb}
\usepackage{amsmath}
\usepackage{xcolor}
\usepackage{cases}
\usepackage{mathtools}
\usepackage{esint}
\usepackage{amsthm}
\usepackage{setspace}
\usepackage{float}
\usepackage{geometry}
\usepackage{verbatim}
\usepackage{mathrsfs}
\usepackage{cite}

\theoremstyle{definition}
\newtheorem{exmp}{Example}[section]

\DeclarePairedDelimiter\floor{\lfloor}{\rfloor}

\newcommand*\diff{\mathop{}\!\mathrm{d}}

\makeatletter
\patchcmd{\@maketitle}{\LARGE \@title}{\fontsize{16}{19.2}\selectfont\@title}{}{}
\makeatother

\usepackage{authblk}

\newsavebox\affbox
\author[1]{\textbf{Adam Johnson}}

\affil[1]{The Fariborz Maseeh Department of Mathematics and Statistics, Portland State Univerisity, 
\newline
Portland, OR, 97207, USA}

\titlespacing\section{0pt}{12pt plus 4pt minus 2pt}{0pt plus 2pt minus 2pt}
\titlespacing\subsection{12pt}{12pt plus 4pt minus 2pt}{0pt plus 2pt minus 2pt}
\titlespacing\subsubsection{12pt}{12pt plus 4pt minus 2pt}{0pt plus 2pt minus 2pt}

\titleformat{\section}{\normalfont\fontsize{10}{15}\bfseries}{\thesection.}{1em}{}
\titleformat{\subsection}{\normalfont\fontsize{10}{15}\bfseries}{\thesubsection.}{1em}{}
\titleformat{\subsubsection}{\normalfont\fontsize{10}{15}\bfseries}{\thesubsubsection.}{1em}{}

\titleformat{\author}{\normalfont\fontsize{10}{15}\bfseries}{\thesection}{1em}{}

\title{\textbf{\huge Solving Diffusion and Wave Equations Meshlessly via Helmholtz Equations}\\}
\date{}    

\begin{document}

\pagestyle{headings}	
\newpage
\setcounter{page}{1}
\renewcommand{\thepage}{\arabic{page}}

\captionsetup[figure]{labelfont={bf},labelformat={default},labelsep=period,name={Figure }}	\captionsetup[table]{labelfont={bf},labelformat={default},labelsep=period,name={Table }}
\setlength{\parskip}{0.5em}
	
\maketitle
	
\noindent\rule{15cm}{0.5pt}
	\begin{abstract}
In this paper, using the approximate particular solutions of Helmholtz equations in \cite{li2009approximation}, we solve the boundary value problems of Helmholtz equations by combining the methods of fundamental solutions (MFS) with the methods of particular solutions (MPS). Then the initial boundary value problems of the time dependent diffusion and wave equations are discretized numerically into a sequence of Helmholtz equations with the appropriate boundary value conditions, which is done by either using the Laplace transform or by using time difference methods. Then Helmholtz problems are solved consequently in an iterative manner, which leads to the solutions of diffusion or wave equations. Several numerical examples are presented to show the efficiency of the proposed methods.

		\let\thefootnote\relax\footnotetext{
			\small $^{*}$ \textit{
				\textit{E-mail address: addy.mathematics@gmail.com}}\\
		                                      }
		\textbf{\textit{Keywords}}: \textit{Radial basis functions; Wave equation; Diffusion equation; Methods of particular solutions, Methods of fundamental solutions }
	\end{abstract}
\noindent\rule{15cm}{0.4pt}

\section{Introduction}

It is well known that the finite element methods (FEM) are widely used in solving differential equations numerically. In using FEM, a domain needs to be meshed geometrically. For instance, in a 2D case, a domain is usually meshed into triangles or rectangles. This process, however, can be complicated and time consuming, especially for irregular domains or three dimensional domains. So in recent years, meshless methods have attracted much attention in engineering, mathematics, and other disciplines for its simplicity in implementation. In other words, instead of generating a mesh, we simply choose points from the domain. The efficiencies of meshless methods have been well documented in the literature (c.f.  \cite{chen2015fast, choi2018meshless, fairweather1998method, golberg1998method, popov2010meshless, moridis1991laplace, zhu1994solving, jankowska2018improved, chen2020novel, zhu2020fictitious}, etc.)

Consider a boundary value problem 

\begin{alignat}{2}
    \mathcal{L} u(\mathbf{x}) &= f(\mathbf{x}), \quad &&\mathbf{x} \in \Omega, \label{1.1.1a} 
   \\
   \mathcal{B} u(\mathbf{x}) &=g (\mathbf{x}), \quad &&\mathbf{x} \in \partial \Omega,  \label{1.1.1b}  
\end{alignat}
where $\mathcal{L}$ is a linear differential operator, $\mathcal{B}$ is a boundary operator, and $\Omega$ is a bounded domain in $\mathbb{R}^s$, $s \geq 2$. Let $\Gamma(\mathbf{x},\mathbf{y})$ be the fundamental solution of $\mathcal{L}$ with singularity at $\mathbf{y}$, i.e.,
$$
\mathcal{L} \Gamma(\mathbf{x},\mathbf{y}) =  \delta(\mathbf{x} - \mathbf{y}), \quad \mathbf{x},\mathbf{y} \in \Omega,
$$
where $\delta$ is the Dirac delta function. To solve problem (\ref{1.1.1a})-(\ref{1.1.1b}) using MFS and MPS, we first need to find a particular solution of (\ref{1.1.1a}). An exact particular solution is not always available, or it is usually given by some singular integral form, which is difficult to calculate for numerical purposes. Therefore, approximate particular solutions are usually desired and commonly used. Suppose that $u_n(\mathbf{x})$ is an approximate particular solution of (\ref{1.1.1a}), namely
$$
\mathcal{L} u_n(\mathbf{x}) \approx f(\mathbf{x}), \quad \mathbf{x} \in \Omega.
$$

Then we consider the following homogeneous problem

\begin{alignat}{2}
    \mathcal{L} u(\mathbf{x}) &= 0, \quad &&\mathbf{x} \in \Omega, \label{1.1.1a homo} 
   \\
   \mathcal{B} u(\mathbf{x}) &=g (\mathbf{x})-u_n(\mathbf{x}), \quad &&\mathbf{x} \in \partial \Omega.  \label{1.1.1b homo}  
\end{alignat}
This homogeneous problem can be solved by using MFS for certain types of differential equations. To be precise, we choose a fictitious domain $\widetilde{\Omega}$ so that $\Omega \subset \widetilde{\Omega}$. Choose $N$ source points $\mathbf{\widetilde{x}_k} \in \partial \widetilde{\Omega}$, where $1 \leq k \leq N$, and then form
\begin{equation}
    u_N(\mathbf{x})=\sum_{k=1}^N c_k \Gamma(\mathbf{x},\mathbf{\widetilde{x}_k}) \label{mfs approx sol}
\end{equation}
for some coefficients {$c_k$} to be determined later. Then $u_N$ satisfies the Equation (\ref{1.1.1a homo}). For $u_N$ to satisfy Equation (\ref{1.1.1b homo}) as close as possible we use N points, called collocation points, $\mathbf{x}_k$, $1 \leq k\leq N$, on $\partial \Omega$, and set
$$
u_N(\mathbf{x}_k)=g(\mathbf{x}_k)-u_n(\mathbf{x}_k)
$$
for $1 \leq k \leq N$, which leads to the following linear system for $\{c_k\}$.

\begin{equation}
\begin{pmatrix}
\Gamma(\mathbf{x}_1,\mathbf{\widetilde{x}}_1) & \Gamma(\mathbf{x}_1,\mathbf{\widetilde{x}}_2) & \cdots & \Gamma(\mathbf{x}_1,\mathbf{\widetilde{x}}_N) \\
\Gamma(\mathbf{x}_2,\mathbf{\widetilde{x}}_1) & \Gamma(\mathbf{x}_2,\mathbf{\widetilde{x}}_2) & \cdots & \Gamma(\mathbf{x}_2,\mathbf{\widetilde{x}}_N) \\
\vdots  & \vdots  & \ddots & \vdots  \\
\Gamma(\mathbf{x}_N,\mathbf{\widetilde{x}}_1) & \Gamma(\mathbf{x}_N,\mathbf{\widetilde{x}}_1) & \cdots & \Gamma(\mathbf{x}_N,\mathbf{\widetilde{x}}_N) 
\end{pmatrix}
\begin{pmatrix}
c_1 \\
c_2 \\
\vdots \\ 
c_N
\end{pmatrix}
=
\begin{pmatrix}
g(\mathbf{x}_1) - u_n(\mathbf{x}_1) \\
g(\mathbf{x}_2)- u_n(\mathbf{x}_2) \\
\vdots \\
g(\mathbf{x}_N) - u_n(\mathbf{x}_N) \\
\end{pmatrix}. \label{mfs system}
\end{equation}

\noindent Once $\{c_k\}$ is determined, $u_N$ in Equation (\ref{mfs approx sol}) is a numerical solution for the boundary value problem (\ref{1.1.1a homo})-(\ref{1.1.1b homo}). Finally, the numerical solution to the original problem (\ref{1.1.1a})-(\ref{1.1.1b}) is given by $u_{N,n}=u_n+u_N$. Such a method of combing MFS with MPS is called the dual reciprocity method (DRM). This method has been widely used to solve boundary value problems of the Laplace, Helmholtz, biharmonic, and many other types of equations (c.f. \cite{fairweather1998method, choi2018meshless, karageorghis2011survey, 2019helm}, etc.).

 In this paper we will use DRM to solve the initial boundary value problems of the diffusion and wave equations along the lines as discussed in the papers \cite{golberg1998method, mulesh}, which the approximate particular solutions of the Helmholtz equations in \cite{li2009approximation} will be used. More precisely, the Laplace transform and time difference methods will be used to discretize diffusion and wave equations.

The organization of this paper is as follows. In \textsection{\ref{helm sec}}, the approximate particular solutions of Helmholtz equations in \cite{li2009approximation} will be described. The Laplace transform will be used in \textsection{\ref{ltd}} to change diffusion and wave equations into Helmholtz problems, and then the inverse Laplace transform will be applied to get the numerical solutions of the original problems. In \textsection{\ref{dtm sec}}, time difference methods will be used to discretize diffusion and wave equation into a sequence of Helmholtz equations which will be solved in an iterative manner. Finally, numerical examples will be presented in \textsection{\ref{exmp sec}} to demonstrate the efficiencies of the methods in this paper.

\section{Approximate Particular Solutions of Helmholtz Equations} \label{helm sec}
A Helmholtz boundary value problem is presented as 
\begin{alignat}{2}
   \Delta u(\mathbf{x}) + \kappa^2 u(\mathbf{x}) &= f(\mathbf{x}), \quad &&\mathbf{x} \in \Omega, \label{helm a}   \\
   u(\mathbf{x}) &= g(\mathbf{x}),  \quad && \mathbf{x} \in \Omega, \label{helm b}
\end{alignat}
\noindent where $\Omega \subset \mathbb{R}^s$, $s \geq 2$, is a bounded domain, $\Delta=\frac{\partial^2}{\partial x_1^2}+ \dots + \frac{\partial^2}{\partial x_s^2}$ is the Laplace differential operator with respect to $\mathbf{x}=(x_1, \dots , x_s)$ in $\mathbb{R}^s$, and $\kappa \in \mathbb{C}$ is a complex constant. The fundamental solution of the Helmholtz equation is given by
\begin{equation}
    \Gamma(\mathbf{x})=\frac{i}{4}  \Big{(}\frac{\kappa}{2 \pi \|\mathbf{x}\|}\Big{)}^{s/2-1} H_{s/2-1}^{(1)}(\kappa \| \mathbf{x} \|), \quad s \geq 2,
\end{equation}
where $H_{s/2-1}^{(1)}$ is a Hankel function of the first kind and $i=\sqrt{-1}$. In particular, when $s=2$ or $s=3$ we have, respectively, 
\begin{equation}
    H_{0}^{(1)}(\mathbf{x})=J_{0}(\mathbf{x})+iY_0(\mathbf{x}), \quad \text{or} \quad H_{1/2}^{(1)}(\mathbf{x})=-i\Big{(} \frac{2}{\pi \mathbf{x}}  \Big{)}^{1/2}e^{i \mathbf{x}},
    \label{hankel}
\end{equation}
where $J_0$ is the Bessel function of the first kind of order $0$, and $Y_0$ is the Bessel function of the second kind of order $0$. Various methods are available in the literature to get approximate particular solutions of Helmholtz equations. A typical one is to use the collocation methods by radial basis functions (RBFs). This method has been used in many other equations as well, such as the Laplace, biharmonic, Stoke's equations, (c.f. \cite{li2010radial, li2009approximation, choi2018meshless, karur1994radial}). Especially if $f(\mathbf{x})$ is a thin plate spline function, the correspondingly exact particular solution of the Helmholtz equation is derived in \cite{mulesh}, and then the collocation methods by using thin plate spline RBFs are used to find approximate particular solutions. Here we describe and use the approximate solutions provided in \cite{li2009approximation}.

In general, an exact solution is given by a singular integral form
\begin{equation}
    u(\mathbf{x}) = \int_{\Omega} \Gamma (\mathbf{x}-\mathbf{y}) f(\mathbf{y}) \diff{y}
\end{equation}
\noindent which, itself as a singular integral, is difficult to evaluate numerically.

Let $\Omega$ be a bounded domain in $\mathbb{R}^s$, $s \geq 2$. Define $\Omega_{\delta}=\Omega + \delta I =\{\mathbf{x}+\mathbf{y}: \mathbf{x} \in \Omega, \mathbf{y}\in \delta I  \}$, where $I=[-1,1]^s$, and $\delta >0$. For any integer $n$, set $I_n(\Omega_{\delta})=\{\mathbf{j} \in \mathbb{Z}^s : [\frac{\mathbf{j}}{n},\frac{\mathbf{j}+\mathbf{1}}{n}]^s \cap \Omega_{\delta} \neq \emptyset   \}$, where $\mathbf{1}=(1,1, \cdots ,1) \in \mathbb{Z}^s$. Choose a radial basis function (RBF) $\phi$ such that
$$
\int_{\mathbb{R}^s} \phi(\mathbf{x}) \diff{\mathbf{x}}=1.
$$

\noindent Then $f(\mathbf{x})$ is approximated in \cite{li2010radial} by using 
$$
f_n(\mathbf{x})=\frac{1}{n^{s(1-\gamma)}} \sum_{\mathbf{j}\in I_n(\Omega_{\delta})} f\Big{(}\frac{\mathbf{j}}{n}\Big{)}\phi(n^{\gamma}\mathbf{x}-\mathbf{j}n^{\gamma-1}), \quad 0 \leq \gamma \leq 1,
$$
\noindent and its order of approximation is derived. For instance, it is shown in \cite{li2010radial} under mild assumptions that if $\gamma=1/2$, then

 $$
 || f_n(\mathbf{x})-f(\mathbf{x}) ||_{C(\Omega)} \leq \frac{c}{\sqrt{n}}||f||_{C^1(\Omega)},
 $$
\noindent where $c$ is independent of $n$ and $f$. Then we consider
\begin{equation}
     \Delta u(\mathbf{x}) + \kappa^2 u(\mathbf{x}) = f_n(\mathbf{x}), \quad \mathbf{x} \in \Omega, \label{helm particular}
\end{equation}
 \noindent whose exact solutions are derived in \cite{li2009approximation}, given by 

\begin{align}
u_n(\mathbf{x})=-\frac{i \pi}{2} \sum_{\mathbf{j}\in I_n(\Omega_{\delta})}f(\frac{\mathbf{j}}{n}) \Bigg{(}  H_0^{(1)} (\kappa  \| \mathbf{x}-\mathbf{j}/n \|)  \int_0^{n^{\gamma} \| \mathbf{x}-\mathbf{j}/n \|} t \phi(t) J_0(\kappa n^{- \gamma }t) \diff{t} \nonumber \\
+ J_0(\kappa  \| \mathbf{x}-\mathbf{j}/n \|) \int_{n^{\gamma} \| \mathbf{x}-\mathbf{j}/n \|}^{\infty}  t \phi(t) H_0^{(1)}(\kappa n^{- \gamma }t) \diff{t}  \Bigg{)}, \label{nonhomohelm approx R2}
\end{align}

\noindent for $\mathbf{x} \in \mathbb{R}^2$ and some constant $\gamma$ (we will use $\gamma =1/2$), and when $\mathbf{x} \in \mathbb{R}^3$,

\begin{align}
u_n(\mathbf{x})= -\frac{i \pi}{2n^3} \sum_{\mathbf{j}\in I_n(\Omega_{\delta})} f(\frac{\mathbf{j}}{n}) n^{-\gamma+ 3 \gamma/2} \Bigg{(}  
\frac{H_{1/2}^{(1)}(\kappa \| \mathbf{x}-\mathbf{j}/n  \|)}{\| \mathbf{x}-\mathbf{j}/n  \|^{1/2}}
\int_0^{n^{\gamma} \| \mathbf{x}-\mathbf{j}/n  \|} t^{3/2} \phi(t) J_{1/2}(\kappa n^{-\gamma}t) \diff{t} \nonumber \\
- \frac{J_{1/2}(\kappa \| \mathbf{x}-\mathbf{j}/n  \|)}{\| \mathbf{x}-\mathbf{j}/n  \|^{1/2}}  \int_{n^{\gamma} \| \mathbf{x}-\mathbf{j}/n  \|}^{\infty}  t^{3/2} \phi(t) H_{1/2}^{(1)}(\kappa n^{-\gamma} t) \diff{t} \Bigg{)}. \label{nonhomohelm approx R3}
\end{align}

\noindent The order of approximation of $u_n(\mathbf{x})$ to the exact solution $u(\mathbf{x})$ is given in \cite{li2009approximation} and roughly described as
$$
|| u_n(\mathbf{x})- u(\mathbf{x}) ||_{C(\Omega)} \leq \frac{c}{\sqrt{n}} || f||_{C^1(\Omega)},
$$
under mild conditions. Then $u_n(\mathbf{x})$ is considered as the approximate particular solution of the original Helmholtz equation and will be used in this paper in the following sections.

In comparison, in using the collocation methods by RBFs for particular solutions, collocation points are chosen from the domain $\Omega$, which leads to a linear system to determine unknown coefficients, exactly as in using the MFS. This generally limits the number of collocation points to be used, since if the number of collocation points is large, it easily leads to an ill-conditioned system and affects the numerical accuracy. Our approximate particular solutions in \cite{li2009approximation} are given by linear summations, and $n$ is allowed to be large as long as the computer software allows, but we certainly like to use small $n$ in concerning the computational time as long as the numerical error is acceptable.

Once $u_n(\mathbf{x})$ is available, we consider the following homogeneous problem.
\begin{alignat}{2}
   \Delta v(\mathbf{x}) + \kappa^2 v(\mathbf{x}) &=0, \quad &&\mathbf{x} \in \Omega, \label{homo-helm a} 
   \\
   v(\mathbf{x}) &=g (\mathbf{x})-u_n(\mathbf{x}), \quad &&\mathbf{x} \in \partial \Omega, \label{homo-helm b}  
\end{alignat}
which is solved by using MFS, as mentioned before. Let $v_p(\mathbf{x})$ be the numerical solution to problem (\ref{homo-helm a})-(\ref{homo-helm b}). Then the solution to the original problem (\ref{helm a})-(\ref{helm b}) is given by
\begin{equation}
u(\mathbf{x})=u_n(\mathbf{x})+v_p(\mathbf{x}). \label{non homo helm approx sol}
\end{equation}

\section{Laplace Transform for Diffusion and Wave Equations} \label{ltd}
Consider an initial boundary value problem (IBVP) of a diffusion equation
\begin{equation}
\frac{1}{k} \frac{\partial u}{\partial t}(\mathbf{x},t) = \Delta u(\mathbf{x},t), \quad \mathbf{x} \in \Omega, \quad t>0, \label{heateq}
\end{equation}
with boundary conditions
\begin{alignat}{2}
    u(\mathbf{x},t) &= g_1(\mathbf{x},t), \quad &&\mathbf{x} \in \partial \Omega_1, \quad t>0, \label{heatbc1}
   \\
   \frac{\partial u}{\partial \mathbf{n}}(\mathbf{x},t) &=g_2 (\mathbf{x},t), \quad &&\mathbf{x} \in \partial \Omega_2,  \quad t>0,  \label{heatbc2}
\end{alignat}
where $\partial \Omega=\partial \Omega_1 \cup \partial \Omega_2$, $\partial \Omega_1 \cap \partial \Omega_2 = \emptyset$, and initial condition
$$
u(\mathbf{x},0)=u_0(\mathbf{x}), \quad \mathbf{x} \in \overline{\Omega}.
$$
The diffusion coefficient $k$ is a constant, and $u_0$, $g_1$, and $g_2$ are known functions.
\noindent Apply the Laplace transform to the above problem (c.f. \cite{golberg1998method}, for instance), where the Laplace transform is defined by
$$
\mathscr{L}\Big{[} u(\mathbf{x},t)\Big{]}=U(\mathbf{x},s)=\int_0^{\infty}u(\mathbf{x},t)e^{-st} \diff{t},
$$
for those $s \in \mathbb{R}$ such that the integral is convergent. Since
$$
\mathscr{L}\Big{[}\frac{\partial u}{\partial t}(\mathbf{x},t)\Big{]}=s\mathscr{L}[u(\mathbf{x},t)]-u(\mathbf{x},0)=s\mathscr{L}[u(\mathbf{x},t)]-u_0(\mathbf{x}),
$$
we get
$$
\frac{1}{k}\Big{(} sU(\mathbf{x},s)-u_0(\mathbf{x}) \Big{)}=\Delta U(\mathbf{x},s)  
$$
i.e.,
\begin{equation}
    \big{(} \Delta-\frac{s}{k} \big{)} U(\mathbf{x},s)=-\frac{u_0(\mathbf{x})}{k}, \quad \mathbf{x} \in \Omega, \label{helm}
\end{equation}
and the corresponding boundary conditions become
\begin{alignat}{2}
   U(\mathbf{x},s) &= G_1(\mathbf{x},s), \quad &&\mathbf{x} \in \partial \Omega_1, \label{helm BC1}
   \\
   \frac{\partial U}{\partial \mathbf{n}}(\mathbf{x},s) &=G_2 (\mathbf{x},s), \quad &&\mathbf{x} \in \partial \Omega_2,  \label{helm BC2}
\end{alignat}
where $G_i(\mathbf{x},s)=\mathscr{L}[g_i(\mathbf{x},t)]$, $i=1,2$. For each $s$, Equation (\ref{helm}) is a Helmholtz equation with differential operator $\Delta-\lambda^2$, $\lambda=\sqrt{s/k}$. If $U(\mathbf{x},s)$ is found, then a solution of the original initial boundary value problem of the diffusion equation is given by
$$
u(\mathbf{x},t)=\mathscr{L}^{-1}\big{[} U(\mathbf{x},s)  \big{]},
$$
where $\mathscr{L}^{-1}$ is the inverse Laplace transform. To evaluate $\mathscr{L}^{-1}[U]$, we use the results in \cite{stehfest1970algorithm} to compute $U(\mathbf{x},s)$ at a finite number of points
\begin{equation}
s_l=\frac{\ln{2}}{t}l, \quad l=1,2,\hdots n_s, \label{helm s value}
\end{equation}
where $t$ is the time at which the solution to our original IBV problem is desired and $n_s$ is some even number. Following the work in \cite{stehfest1970algorithm}, the numerical approximation $\widetilde{u}(\mathbf{x},t)$ to $u(\mathbf{x},t)$ is given by
\begin{equation}
    \widetilde{u}(\mathbf{x},s)=\frac{\ln{2}}{t} \sum_{l=1}^{n_s} \alpha_l U(\mathbf{x},s_l), \label{inverse laplace formula}
\end{equation}
\noindent where
\begin{equation}
\alpha_l=(-1)^{n_s/2+l} \sum_{i=\floor{\frac{l+1}{2}}}^{\min{\{l, n_s/2}\}} \frac{i^{n_s/2}(2i)!}{i!(n_s/2-i)!(i-1)!(l-i)!(2i-l)!}. \label{inv laplace constant}
\end{equation}
The above formula depends on the choices of $n_s$ and it is typical to choose $n_s=10$ or $18$.

The method can be similarly applied to wave equations. Consider the IBVP of a wave equation given by

\begin{equation}
\begin{aligned}
    \frac{1}{c^2} \frac{\partial^2 u}{\partial t^2}(\mathbf{x},t) &= \Delta u(\mathbf{x},t), && \mathbf{x} \in \Omega, \quad & t > 0,  \\
    \mathcal{B}(u(\mathbf{x},t)) &= g(\mathbf{x},t), && \mathbf{x} \in \partial \Omega, \quad & t > 0, \\
    \frac{\partial u}{\partial t}(\mathbf{x},0) &= v_0(\mathbf{x}), && \mathbf{x} \in \Omega, \quad & t = 0, \\
    u(\mathbf{x},0) &= u_0(\mathbf{x}), && \mathbf{x} \in \Omega, \quad & t = 0.
\end{aligned}
\end{equation}

  \noindent where $\Omega \subset \mathbb{R}^s$, $s \geq 2$, is a bounded domain, $\mathcal{B}$ is a boundary operator, and $c$ is the the propagation velocity.

Applying the Laplace transform to the above problem, we arrive at
\begin{alignat}{2}
    \Big{(}  \frac{s}{c} \Big{)}^2 U(\mathbf{x},s) - s u_0(\mathbf{x}) - v_0(\mathbf{x}) &=\Delta U(\mathbf{x},s), \quad &&\mathbf{x} \in \Omega,   \label{LT heat eq v1}  \\
    \mathcal{B}(U(\mathbf{x},s))&=G(\mathbf{x},s), \quad  && \mathbf{x} \in \partial \Omega,
  \end{alignat}
  \noindent where $G(\mathbf{x},s) = \mathscr{L}[g(\mathbf{x},t)]$. Define $F(\mathbf{x},s)= s u_0(\mathbf{x}) + v_0(\mathbf{x})$, and $\lambda=s/c$, then we can rewrite Equation (\ref{LT heat eq v1}) as 
  \begin{equation*}
       \Delta U(\mathbf{x},s)-\lambda^2 U(\mathbf{x},s) = F(\mathbf{x},s), \quad \mathbf{x} \in \Omega.
  \end{equation*}
  
  \noindent For each fixed $s
  $, we again get a boundary value problem of a Helmholtz equation. And if $U(\mathbf{x},s)$ is available, we use the inverse Laplace transform of $U(\mathbf{x},s)$ to get the solution as discussed above.

\section{Difference Method for the Diffusion  and Wave Equations} \label{dtm sec}
Here we discuss difference in time methods for the diffusion IBVP (\ref{heateq})-(\ref{heatbc2}) with $\Omega = \Omega_1$, i.e., we only consider the Dirichlet boundary condition (\ref{heatbc1}). Let $\tau > 0$ and define $u_n(\mathbf{x})=u(\mathbf{x}, n \tau)$, $n \geq 0$, and approximate the derivative $u_t$ by
  \begin{equation}
      u_t(\mathbf{x}, n \tau) \approx \frac{u(\mathbf{x}, n \tau)-u(\mathbf{x}, (n-1) \tau)}{\tau}. \label{ut approx}
  \end{equation}
  Use an approximation $\widetilde{u}_n(\mathbf{x})$ to $u_n(\mathbf{x})$ by requesting 
 \begin{alignat}{2}
    \Delta \widetilde{u}_n(\mathbf{x}) - k \frac{\widetilde{u}_n(\mathbf{x})-\widetilde{u}_{n-1}(\mathbf{x})}{\tau} &= f(\widetilde{u}_{n-1}(\mathbf{x})), \quad && \mathbf{x} \in \Omega, \quad n \geq 1, \label{approx fdm heat a}\\
    \widetilde{u}_0(\mathbf{x})&=u_0(\mathbf{x}), \quad  && \mathbf{x} \in \Omega, \label{approx fdm heat b}  \\ 
    \widetilde{u}_n(\mathbf{x})&=g(\mathbf{x}), \quad &&\mathbf{x} \in \partial \Omega. \label{approx fdm heat c}
  \end{alignat}
 
 \noindent Rewrite Equation (\ref{approx fdm heat a}) as
  \begin{equation}
      \Delta \widetilde{u}_n(\mathbf{x}) - \frac{k}{\tau} \widetilde{u}_n(\mathbf{x}) = - \frac{k}{\tau} \widetilde{u}_{n-1}(\mathbf{x})+f(\widetilde{u}_{n-1}(\mathbf{x})), \quad n \geq 1 .
  \end{equation}
  Set $\lambda = k/ \tau$, $k>0$, and $h_n(\mathbf{x})= - \frac{k}{\tau} \widetilde{u}_{n-1}(\mathbf{x})+f(\widetilde{u}_{n-1}(\mathbf{x}))$. Then we have the inhomogeneous Helmholtz equation given by
  \begin{equation}
      \Delta \widetilde{u}_n(\mathbf{x}) - \lambda^2 \widetilde{u}_n(\mathbf{x}) = h_n(\mathbf{x}), \quad \mathbf{x} \in \Omega, \label{fdm helm eq}
  \end{equation}
  along with the initial condition (\ref{approx fdm heat b}) and boundary condition (\ref{approx fdm heat c}), which can be solved as discussed before.

The method can be similarly applied to the wave IBVP with a Dirichlet boundary condition. With a fixed $\tau > 0$, and $t_n=n\tau$ for $n \geq 1$, we let $u_n=u(\mathbf{x},t_n)$ for the sake of notation. Then for $n \geq 2$,
  $$
  \frac{\partial^2 u}{\partial t^2}(\mathbf{x},t_n) \approx \frac{u_{n}-2u_{n-1}+u_{n-2}}{\tau^2}.
  $$
 The wave equation can now be approximated as
  $$
\Delta u_n = \frac{1}{c^2} \frac{u_n-2u_{n-1}+u_{n-2}}{\tau^2}, \quad n \geq 2,
  $$
  i.e.,
  \begin{equation}
      \Delta u_{n} - \frac{1}{c^2 \tau^2} u_n = \frac{1}{c^2 \tau^2} \big{(} u_{n-2} -2 u_{n-1} \big{)}, \quad n \geq 2. \label{time diff wave eq}
  \end{equation}
  Now, using the approximation
  $$
  u_t(\mathbf{x},t_n) \approx \frac{u_{n+1}-u_n}{\tau}
  $$
  we get the following initial conditions
  \begin{equation}
      u_0(\mathbf{x})=u(\mathbf{x},0), \quad u_1(\mathbf{x})=u_0(\mathbf{x})+v_0(\mathbf{x}) \tau, \quad \mathbf{x} \in \Omega, \label{time diff wave ic}
  \end{equation}
  and the boundary condition becomes
  \begin{equation}
      u_{n+1}(\mathbf{x})=g_{n+1}(\mathbf{x}), \quad \mathbf{x} \in \partial \Omega. \label{time diff wave bc}
  \end{equation}
  Notice problem (\ref{time diff wave eq})-(\ref{time diff wave bc}) is a non-homogeneous Helmholtz boundary value problem and thus it can be solved by the methods mentioned above.

\section{Numerical Examples} \label{exmp sec}

In this section we present a few examples that demonstrate the efficiencies of the numerical methods described in the above sections. First, we begin with an example of solving a BVP of a Helmholtz equation by using MFS in a 2D case.

\begin{exmp} \label{helm1ex} Consider a homogeneous Helmholtz problem
\begin{alignat}{2}
    \Delta u(x,y)-3u(x,y) &= 0, \quad &&(x,y) \in \Omega,  \label{helm ex 1a}    \\
    u(x,y)&=e^{-2x} \sin{y}, \quad  &&(x,y) \in \partial \Omega, \label{helm ex 1b}
  \end{alignat} 

\noindent where $\Omega$ is the domain whose boundary is given by the polar equation
$$
r=1-\frac{1}{3} \cos{(4 \theta)},
$$
\noindent where $\theta \in [0,2 \pi]$. The exact solution to this problem is given by $u_e(x,y)=e^{-2x} \sin{y}$. We will choose $N$ collocation points on $\partial \Omega$ by plugging in $N$ equally spaced points of $[0,2 \pi]$ into the above polar equation. We will choose $N$ source points to be equally spaced points on the circle of radius $a$ centered at the origin (see Figure \ref{helm1fig}).

\begin{figure}[H]
  \centering
    \includegraphics[width=.6\textwidth]{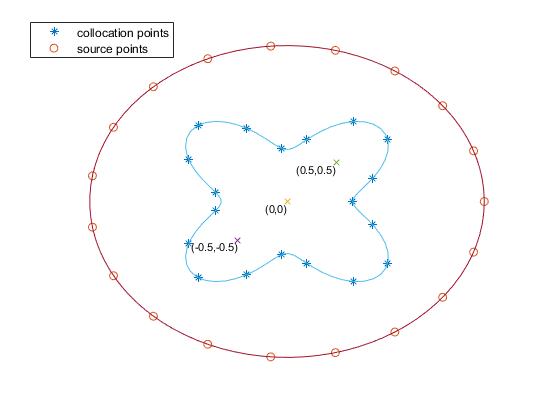}
    \caption{Collocation, interior, and source points for Example \ref{helm1ex}.}
    \label{helm1fig}
\end{figure}

\noindent Then setting $u_n=0$ in Equation (\ref{non homo helm approx sol}) gives us the numerical solution to problem (\ref{helm ex 1a})-(\ref{helm ex 1b}). Equation (\ref{helm ex 1a}) is a modified Helmholtz equation for which the maximum principle applies. So to estimate the numerical error we choose $150$ points $\mathbf{z}_k \in \partial \Omega$ corresponding to equally spaced points in $[0,2 \pi]$ and calculate
$$
\text{error} =\max_{1 \leq k \leq 150} \vert u_e(\mathbf{z}_k) -u_N(\mathbf{z}_k) \vert.
$$
\noindent Furthermore, we compare the numerical solution with the exact solution at the interior points of $\Omega$ shown in Figure \ref{helm1fig}: $(0,0),(-0.5,-0.5)$, and $(0.5,0.5)$. Then our results are presented below for different values of $N$ and $a$.

 \begin{table}[H]
 \centering
 \caption{Numerical errors on $\partial \Omega$ for Example \ref{helm1ex}}
\begin{tabular}{ |p{3cm}||p{3cm}|p{3cm}|p{3cm}|  } 
 \hline
 $a$ Value& $N=30$ & $N=40$ &$N=50$\\
 \hline
 $a=5$   & 5.7834e-08  &  1.7858e-04 &   1.8384e-05\\
 $a=6$&   6.4251e-08  &  1.5524e-04 & 3.0641e-04\\
 $a=7$ & 7.6068e-08 & 3.2974e-04&  2.1189e-04\\
 \hline   
\end{tabular} 
\end{table}

 \begin{table}[H]
 \centering
 \caption{Numerical errors at $(0,0) \in \Omega$ for Example \ref{helm1ex}}
\begin{tabular}{ |p{3cm}||p{3cm}|p{3cm}|p{3cm}|  } 
 \hline
 $a$ Value& $N=30$ & $N=40$ &$N=50$\\
 \hline
 $a=5$   & 1.2697e-09 & 3.7395e-05 &   2.0255e-06\\
 $a=6$  & 3.6321e-09  &  2.6167e-05 & 4.5180e-06\\
 $a=7$ &  7.2834e-09&  3.8662e-05&  7.6944e-06 \\
 \hline   
\end{tabular} 
\end{table}

 \begin{table}[H]
 \centering
 \caption{Numerical errors at $(-0.5,-0.5) \in \Omega$ for Example \ref{helm1ex}}
\begin{tabular}{ |p{3cm}||p{3cm}|p{3cm}|p{3cm}|  } 
 \hline
 $a$ Value& $N=30$ & $N=40$ &$N=50$\\
 \hline
 $a=5$   & 3.3531e-09 &  1.4115e-05 &   3.6029e-06\\
 $a=6$ &  1.0839e-08 & 5.9511e-05 & 3.6455e-05\\
 $a=7$ & 2.6621e-09 & 8.7411e-05&  6.5502e-06 \\
 \hline   
\end{tabular} 
\end{table}

 \begin{table}[H]
 \centering
 \caption{Numerical errors at $(0.5,0.5) \in \Omega$ for Example \ref{helm1ex}}
\begin{tabular}{ |p{3cm}||p{3cm}|p{3cm}|p{3cm}|  } 
 \hline
 $a$ Value& $N=30$ & $N=40$ &$N=50$\\
 \hline
 $a=5$   & 1.8178e-10 & 7.0093e-05 &   5.2808e-06\\
 $a=6$&   1.2718e-09  & 1.5910e-05 &  9.0011e-05\\
 $a=7$ & 1.1732e-08 & 3.2128e-06 &  6.7437e-05 \\
 \hline   
\end{tabular} 
\end{table}

\begin{figure}[H]
  \begin{minipage}[b]{0.5\textwidth}
    \includegraphics[width=\textwidth]{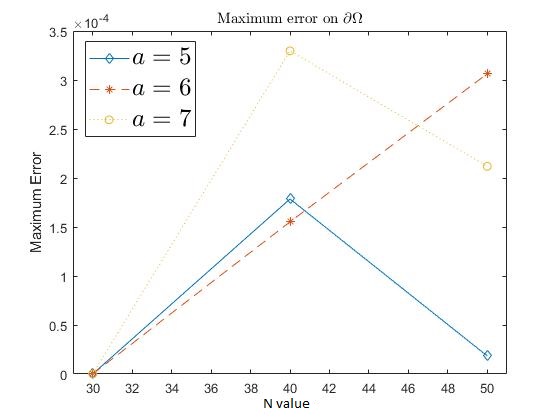}
  \end{minipage}
  \hfill
  \begin{minipage}[b]{0.5\textwidth}
    \includegraphics[width=\textwidth]{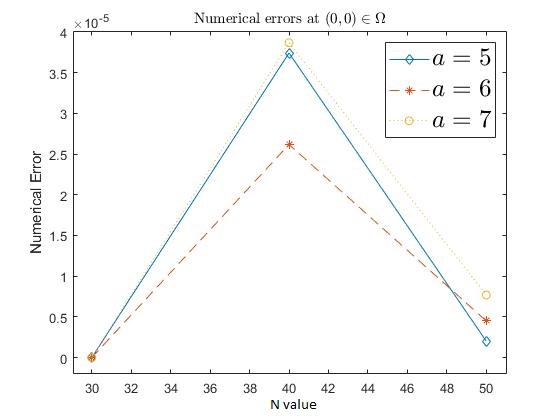}
  \end{minipage}
\end{figure}

\begin{figure}[H]
  \begin{minipage}[b]{0.5\textwidth}
    \includegraphics[width=\textwidth]{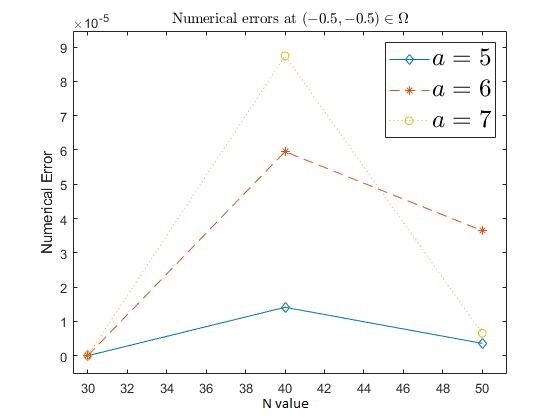}
  \end{minipage}
  \hfill
  \begin{minipage}[b]{0.5\textwidth}
    \includegraphics[width=\textwidth]{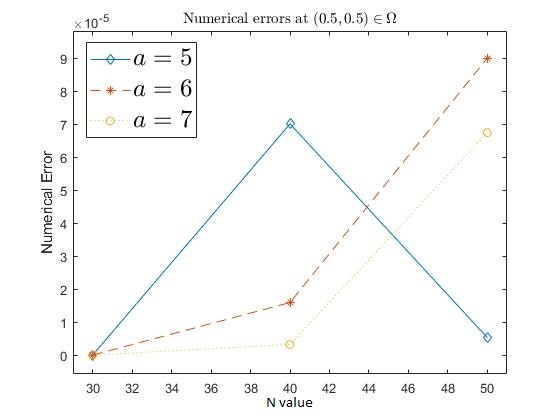}
  \end{minipage}
  \caption{Example \ref{helm1ex} errors on $\partial \Omega$ and $3$ interior points}
\end{figure}

\end{exmp}

Next we use MFS and MPS to solve a non-homogeneous Helmholtz problem.

\begin{exmp} \label{nonhomo helm ex1}
Consider
\begin{alignat}{2}
    \Delta u(x,y)-u(x,y) &= 2xe^y, \quad &&(x,y) \in \Omega,  \label{nonhomohelm ex 1a}    \\
    u(x,y)&=xye^y, \quad  &&(x,y) \in \partial \Omega. \label{nonhomohelm ex 1b}
  \end{alignat} 
where $\Omega = \{(x,y): x^2+y^2 < 1\}$. The exact solution to this problem is given by $u_e(x,y)=xye^y$. For the non-homogeneous Helmholtz equation, we use the approximate particular solutions discussed in \textsection{\ref{helm sec}}. Here we use different RBFs, $\phi_1(r)=\sqrt{\frac{c}{\pi}}e^{-cr^2}$ and $\phi_2(r)=\frac{35}{32}(1-r^2)^3 \chi_{[-1,1]}(r)$, where $\chi$ is the characteristic function. And use points of the form $\frac{\mathbf{j}}{n}=(\frac{k}{n},\frac{m}{n})$ in $\Omega_{\delta}$, $k,m \in \mathbb{Z}$ to get the approximate solutions from Equation (\ref{nonhomohelm approx R2}).

Next, we use the methods in MFS to solve problem (\ref{homo-helm a})-(\ref{homo-helm b}), as discussed before. We obtain $N$ collocation points by plugging in $N$ equally spaced points of $[0, 2 \pi]$ into $(\cos{t},\sin{t})$, and use a fictitious domain $\partial{\widetilde{\Omega}}$ whose boundary is given by the polar equation 
$$
r=5-\cos{6 \theta},
$$
for $\theta \in [0, 2 \pi]$. We get $N$ source points by plugging into the above polar equation. 

\begin{figure}[H]
\centering
  \includegraphics[width=.6\textwidth]{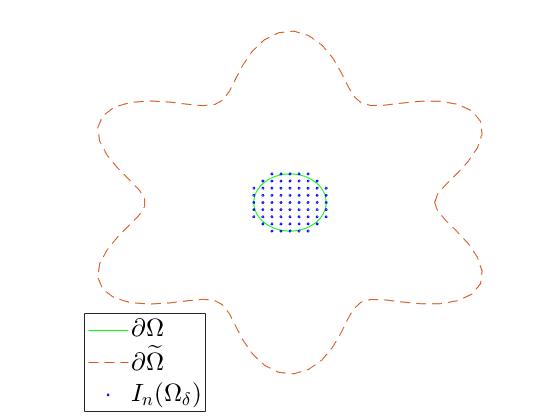} 
  \caption{Domains and points of $I_n(\Omega_{\delta})$ for Example \ref{nonhomo helm ex1}}
\end{figure}

Then the numerical solution to the original problem is given by Equation (\ref{non homo helm approx sol}). To estimate the error we calculate
$$
\text{error}=\max_{1 \leq i \leq 450}{\vert u_e(\mathbf{z}_i)-u_n(\mathbf{z}_i) \vert}
$$
where $u_n$ is the numerical solution to the original problem, and $\mathbf{z}_i$, $1 \leq i \leq 450$, consists of $350$ points in $\Omega$ and $100$ points on $\partial \Omega$. The interior points for estimating the errors are obtained by plugging into the polar equation $r=(1/14)\theta$, for $2 \pi \leq \theta \leq 4 \pi$, and the boundary points are obtained by plugging into $(\cos{t},\sin{t})$. The results are given in the tables below.

 \begin{table}[H] 
 \centering
 \caption{Errors for Example \ref{nonhomo helm ex1} with $\phi(r)=\sqrt{\frac{c}{\pi}}e^{-cr^2}$, $n=14$, and $\delta=0.2$}
\begin{tabular}{ |p{3cm}||p{3cm}|p{3cm}|p{3cm}|  } 
 \hline
 $c$ Value& $N=30$ & $N=40 $ &$N=50$\\
 \hline

 $c=3.77$ & 9.8458e-04 &9.0813e-04 & 9.4048e-04 \\ 
 
   $c=3.8$  & 6.8730e-04 & 6.2337e-04
 &6.6379e-04  \\

 $c=3.82$ &4.9107e-04 & 4.3522e-04 & 4.8256e-04\\

 \hline   
\end{tabular}  
\end{table}

 \begin{table}[H] 
 \centering
 \caption{Errors for Example \ref{nonhomo helm ex1} with $\phi(r)=\sqrt{\frac{c}{\pi}}e^{-cr^2}$, $n=24$, and $\delta=0.1$}
\begin{tabular}{ |p{3cm}||p{3cm}|p{3cm}|p{3cm}|  } 
 \hline
 $c$ Value& $N=30$ & $N=40 $ &$N=50$\\
 \hline

 $c=4$ &1.9549e-03 & 1.9259e-03 & 1.9475e-03 \\ 
 
   $c=4.2$  & 9.4343e-04
 & 8.7272e-04& 9.5018e-04
\\

 $c=4.4$ &2.5448e-03  & 2.4606e-03 &2.5349e-03  \\

 \hline   
\end{tabular}  
\end{table}

 \begin{table}[H] 
 \centering
 \caption{Errors for Example \ref{nonhomo helm ex1} with $\phi_2(r)=\frac{35}{32}(1-r^2)^3 \chi_{[-1,1]}(r)$ and $\delta=0.2$}
\begin{tabular}{ |p{3cm}||p{3cm}|p{3cm}|p{3cm}|  } 
 \hline
 $n$ Value& $N=30$ & $N=40 $ &$N=50$\\
 \hline

 $n=16$ & 1.4661e-03  & 1.4744e-03  &1.4839e-03  \\ 
 
   $n=30$  & 2.3721e-03
 &2.3513e-03   & 2.3640e-03 
\\

 $n=50$ & 9.7023e-04 & 9.7664e-04& 1.0239e-03  \\

 \hline   
\end{tabular}  
\end{table}

\begin{figure}[H]
  \begin{minipage}[b]{0.5\textwidth}
    \includegraphics[width=\textwidth]{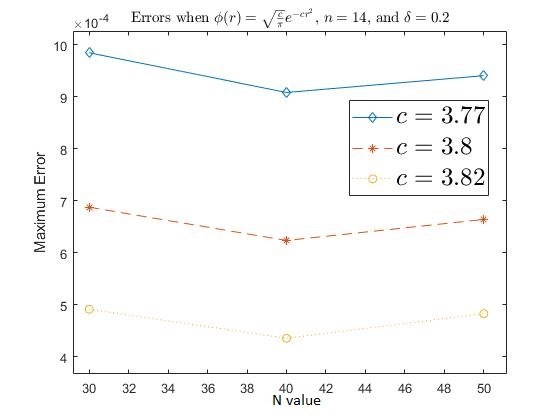}
  \end{minipage}
  \hfill
  \begin{minipage}[b]{0.5\textwidth}
    \includegraphics[width=\textwidth]{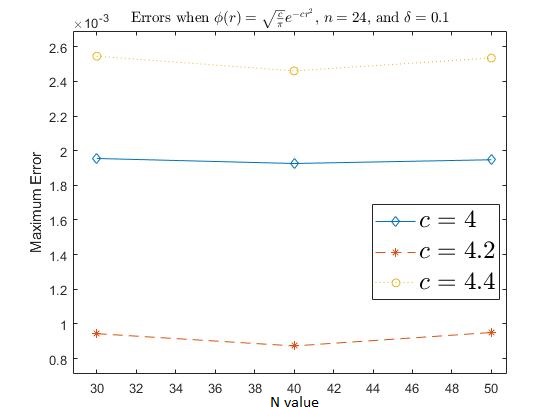}
  \end{minipage}
\end{figure}

\begin{figure}[H]
\centering
  \includegraphics[width=.5\textwidth]{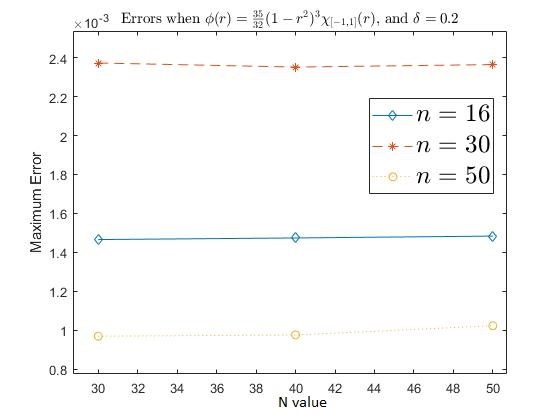} 
  \caption{Example \ref{nonhomo helm ex1} errors for different RBFs, $n$, and $c$ values}
\end{figure}

\end{exmp}

Now we use the particular solution in Equation (\ref{nonhomohelm approx R3}) for a 3D problem.

\begin{exmp} \label{nonhomohelm R3 ex}
Consider the following problem in $\mathbb{R}^3$

\begin{alignat}{2}
    \Delta u(x,y,z)-u(x,y,z) &= 2ye^z, \quad &&(x,y,z) \in \Omega,  \label{nonhomohelm r3 a}    \\
    u(x,y,z)&=x^2ye^z, \quad  &&(x,y,z) \in \partial \Omega. \label{nonhomohelm r3 b}
  \end{alignat} 
\noindent where $\Omega=\{(x,y,z): x^2+y^2+z^2 <1 \}$ is the unit sphere. The exact solution of this problem is $u_e(x,y,z)=x^2ye^z$. We use the Gaussian function $\phi(r)=\sqrt{0.1/\pi}e^{-0.1r^2}$ as our RBF. Then we use points of the form $\frac{\mathbf{j}}{n}=(\frac{k}{n},\frac{m}{n}, \frac{l}{n})$ in $\Omega_{\delta}$, $k,m,l \in \mathbb{Z}$ to get the approximate solutions from Equation (\ref{nonhomohelm approx R3}). We use $N$ collocation points and $N$ source points on a sphere of radius $4$ centered at the origin. The the numerical solution to problem (\ref{nonhomohelm r3 a})-(\ref{nonhomohelm r3 b}) is given by Equation (\ref{non homo helm approx sol}).

\begin{figure}[H]
\centering
  \includegraphics[width=.8\textwidth]{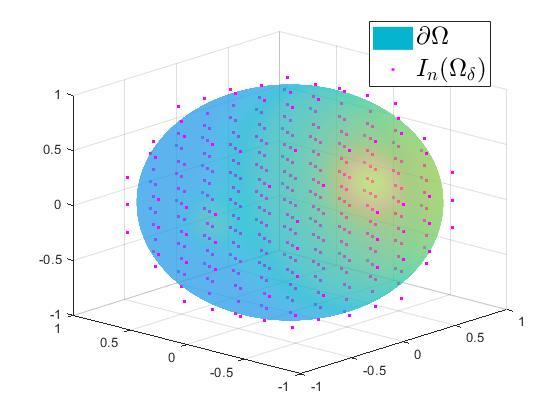} 
  \caption{The boundary $\partial \Omega$ and points in $I_n(\Omega_{\delta})$ for Example \ref{nonhomohelm R3 ex}}
\end{figure}

Error estimates are calculated as in the previous examples. We calculate the maximum error at $686$ boundary points and the numerical error at $3$ interior points. The results are given below.

 \begin{table}[H] 
 \centering
 \caption{Errors on $\partial \Omega$ for Example \ref{nonhomohelm R3 ex} with $\delta=0.2$}
\begin{tabular}{ |p{3cm}||p{3cm}|p{3cm}|p{3cm}|  } 
 \hline
 $n$ Value& $N=155$ & $N=176 $ &$N=203$\\
 \hline
 $n=5$ & 8.8797e-04 &5.8049e-04 &1.0279e-03 \\ 
   $n=10$ & 5.7091e-04& 3.0869e-04  &8.3734e-05 \\
 $n=20$ &4.8601e-04 &2.5440e-04 & 1.0339e-04 \\ 
 \hline   
\end{tabular} 
\end{table}

 \begin{table}[H] 
 \centering
 \caption{Errors at $(0.26,0,-0.15) \in \Omega$ for Example \ref{nonhomohelm R3 ex} with $\delta=0.2$}
\begin{tabular}{ |p{3cm}||p{3cm}|p{3cm}|p{3cm}|  } 
 \hline
 $n$ Value& $N=155$ & $N=176 $ &$N=203$\\
 \hline
 $n=5$ &1.3216e-04 &2.9446e-05  &1.2316e-04 \\ 
   $n=10$ &7.63554e-05 & 1.2761e-06 &1.3956e-05\\
 $n=20$ &3.4145e-05 &1.5425e-05 & 1.3176e-05  \\ 
 \hline   
\end{tabular} 
\end{table} 

\begin{table}[H] 
 \centering
 \caption{Errors at $(0.26,0,0.15) \in \Omega$ for Example \ref{nonhomohelm R3 ex} with $\delta=0.2$}
\begin{tabular}{ |p{3cm}||p{3cm}|p{3cm}|p{3cm}|  } 
 \hline
 $n$ Value& $N=155$ & $N=176 $ &$N=203$\\
 \hline
 $n=5$ &1.5254e-04 & 1.9090e-05 & 1.2736e-04\\ 
   $n=10$ &1.0593e-04  & 8.1219e-06 & 2.2554e-05\\
 $n=20$ &6.2453e-05 &2.6541e-05 & 2.2394e-05 \\ 
 \hline   
\end{tabular} 
\end{table}

\begin{table}[H] 
 \centering
 \caption{Errors at $(0,0,0.3) \in \Omega$ for Example \ref{nonhomohelm R3 ex} with $\delta=0.2$}
\begin{tabular}{ |p{3cm}||p{3cm}|p{3cm}|p{3cm}|  } 
 \hline
 $n$ Value& $N=155$ & $N=176 $ &$N=203$\\
 \hline
 $n=5$ & 8.2412e-05 & 4.6458e-05 & 5.4867e-05 \\ 
   $n=10$ &6.4257e-05 & 2.32639e-05 &1.3766e-05\\
 $n=20$ &5.10585e-05 &2.5940e-05 & 1.3930e-05 \\ 
 \hline   
\end{tabular} 
\end{table}

\begin{figure}[H]
  \begin{minipage}[b]{0.5\textwidth}
    \includegraphics[width=\textwidth]{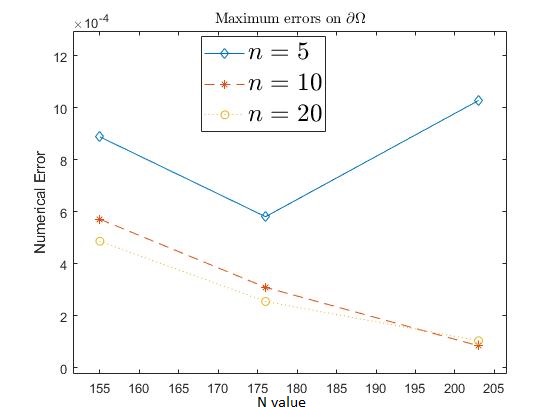}
  \end{minipage}
  \hfill
  \begin{minipage}[b]{0.5\textwidth}
    \includegraphics[width=\textwidth]{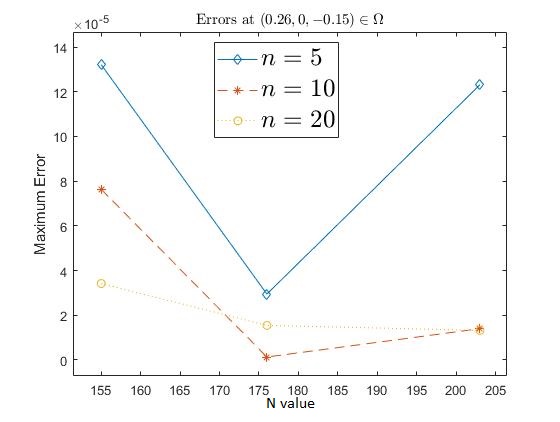}
  \end{minipage}
\end{figure}

\begin{figure}[H]
  \begin{minipage}[b]{0.5\textwidth}
    \includegraphics[width=\textwidth]{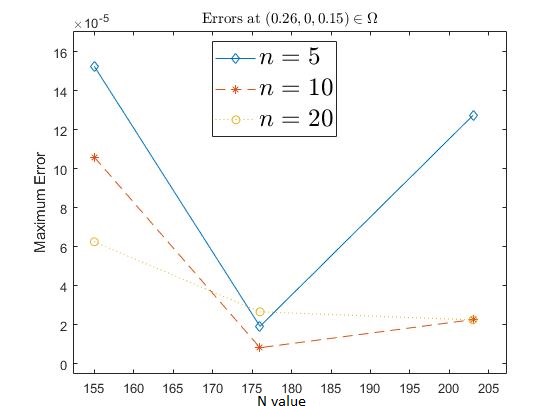}
  \end{minipage}
  \hfill
  \begin{minipage}[b]{0.5\textwidth}
    \includegraphics[width=\textwidth]{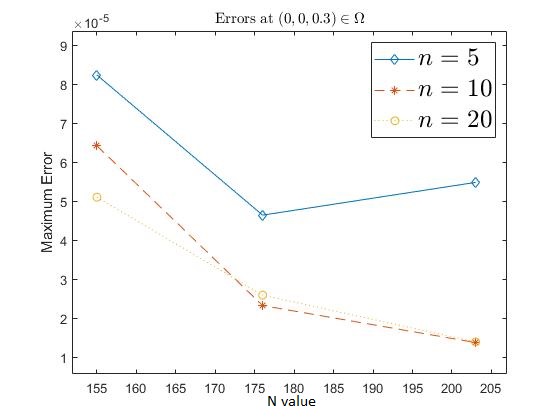}
  \end{minipage}
  \caption{Example \ref{nonhomohelm R3 ex} errors on $\partial \Omega$ and $3$ interior points}
\end{figure}

\end{exmp}

Next we use the Laplace transform for a diffusion problem.

\begin{exmp} \label{heatex}
Consider the following diffusion problem on a square with boundary temperatures at zero. A similar example can be found in \cite{chenlecture}.

\begin{equation}
\begin{aligned}
    \frac{1}{k} \frac{\partial u}{\partial t}(x,y,t) &= \Delta u(x,y), && (x,y) \in \Omega, \quad & t > 0, \label{heat ex 1 eq} \\ 
    u(x,y,t) &= 0, && (x,y) \in \partial \Omega, \quad & t > 0, 
    \\
    u(x,y,0) &= 1, && (x,y) \in \Omega, \quad & t = 0. 
\end{aligned} 
\end{equation}

  \noindent where $\Omega=(-0.1,0.1) \times (-0.1,0.1)$ and $k=5.8 \times 10^{-7}$ $m^2/s$. The analytical solution to this problem is given by
  \begin{equation}
  u_e(x,y,t)=\frac{16}{\pi^2} \sum_{n,m=1}^{\infty} A_{n,m} e^{\lambda_{n,m}t} \cos{\frac{(2n+1)\pi x}{0.2}} \cos{\frac{(2m+1)\pi y}{0.2}}, \label{heatex anal sln}
  \end{equation}
  where
  $$
  A_{n,m}=\frac{(-1)^{n+m}}{(2n+1)(2m+1)},
  $$
  and
  $$
  \lambda_{n,m}=-\frac{k \pi^2}{0.04} \big{(} (2n+1)^2 +(2m+1)^2 \big{)}.
  $$
  Since the initial condition is harmonic we can apply the change of variable $v(x,y,t)=u(x,y,t)-u(x,y,0)$, which leads to the following problem
\begin{equation}
\begin{aligned}
    \frac{\partial v}{\partial t}(x,y,t) &= \Delta v(x,y), && (x,y) \in \Omega, \quad & t > 0, \\
    v(x,y,t) &= -1, && (x,y) \in \partial \Omega, \quad & t > 0, \\
    v(x,y,0) &= 0, && (x,y) \in \Omega, \quad & t = 0.
\end{aligned} \label{heat ex 1 changed eq}
\end{equation}

  \noindent Now we take the Laplace transform of the above problem, as mentioned in \textsection{\ref{ltd}}, to obtain the following Helmholtz problem
  
  \begin{alignat}{2}
    \Delta V(x,y,s_l)-s_l V(x,y,s_l) &= 0, \quad &&(x,y) \in \Omega,  \label{LT heat ex 1}   \\
    V(x,y,s_l)&=-1/s_l, \quad  &&(x,y) \in \partial \Omega, \label{LT heat ex 1 bc} 
  \end{alignat}
  where $V(x,y,s_l)=\mathscr{L}[v(x,y,t)]$, and the values for $s_l$, $1 \leq l \leq n_s$, are given in Equation (\ref{helm s value}) where we choose $n_s=10$ and $18$, and a time $t=9000$ to evaluate the original problem at. Then for each $s_l$ we solve problem (\ref{LT heat ex 1})-(\ref{LT heat ex 1 bc}) using the MFS and MPS where we use $40$ equally spaced points on $\partial \Omega$ as the collocation points and $40$ equally spaced source points on the circle centered at the origin with radius $2$ (See Figure \ref{heatexpic}). Next we approximate the inverse Laplace transform by Equation (\ref{inverse laplace formula}) to get the numerical solution $\widetilde{v}(x,y,t)$ to the problem (\ref{heat ex 1 changed eq}). Finally the numerical solution to the original problem (\ref{heat ex 1 eq}) is given by $\widetilde{u}(x,y,t)=\widetilde{v}(x,y,t)+u(x,y,0)$. 
  
  We compare the exact solution with the numerical solution at $6$ points $(x,y) \in \Omega$, and calculate the error as follows
  $$
  \text{Error}= |u(x,y,9000)-\widetilde{u}(x,y,9000) |.
  $$
  
  The results are given in the tables below.
  
 \begin{table}[H] 
 \centering
 \caption{True and numerical solutions to Example \ref{heatex} with $n_s=10$}
\begin{tabular}{ |p{3cm}||p{3cm}|p{3cm}|p{3cm}|  } 
 \hline
 $(x,y)\in \Omega$ & True \newline Solution & Numerical \newline Solution & Error  \\
 \hline
 $(-0.01,0.07)$ &   0.7014      & 0.7521 & 0.0507\\ 
   $(-0.01,0.04)$ & 0.7682  &  0.8034
 & 0.0351\\
 $(-0.01,0.01)$ &  0.8008  &   0.8026  & 1.7366e-03  \\ 
  $(-0.01,-0.01)$ &   0.8008    &  0.8001    &  7.3928e-04 \\
   $(-0.01,-0.04)$ &   0.7682  &   0.7966   &  0.0283\\
    $(-0.01,-0.07)$ &0.7014   &  0.7929  &  0.0915\\
  
 \hline   
\end{tabular} 
\end{table}

  \begin{table}[H] 
 \centering
 \caption{True and numerical solutions to Example \ref{heatex} with $n_s=18$}
\begin{tabular}{ |p{3cm}||p{3cm}|p{3cm}|p{3cm}|  } 
 \hline
 $(x,y)\in \Omega$ & True \newline Solution & Numerical \newline Solution & Error  \\
 \hline
 $(-0.01,0.07)$ &   0.7014      & 0.7921 & 0.0906\\ 
   $(-0.01,0.04)$ & 0.7682  &  0.7934
 & 0.0251\\
 $(-0.01,0.01)$ &  0.8008  &    0.7926 &  8.2685e-03 \\ 
  $(-0.01,-0.01)$ &   0.8008    & 0.7901   & 0.0107   \\
   $(-0.01,-0.04)$ &   0.7682   &    0.7866  & 0.0183  \\
    $(-0.01,-0.07)$ &0.7014   &  0.7829 & 0.0814 \\
  
 \hline   
\end{tabular} 
\end{table}

\begin{figure}[H]
\centering
  \includegraphics[width=.5\textwidth]{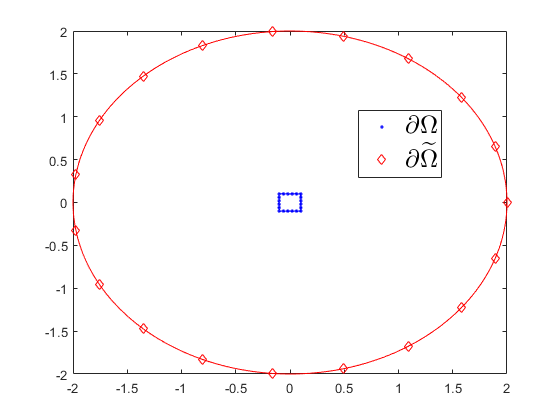} 
  \caption{Collocation and source points for Example \ref{heatex}} \label{heatexpic}
\end{figure}
  
  \end{exmp}

\begin{exmp} \label{heatfdm}
  We now revisit problem (\ref{heat ex 1 eq}) in Example \ref{heatex}, but use the difference in time methods. Recall the analytical solution to this problem is given by Equation (\ref{heatex anal sln}).
  
  First we choose $M$ equally spaced points on the time interval $[0,T]$ so that we can form the approximation of the time derivative $\frac{\partial u}{\partial t}(x,y,n \tau)$ given in Equation (\ref{ut approx}), where $\tau=t/n$, and $n=1,2,\cdots ,M$. Next, we rearrange terms so that our problem becomes the inhomogeneous Helmholtz problem mentioned in Equation (\ref{fdm helm eq}). Then we use the methods in \textsection{\ref{helm sec}} to solve the Helmholtz equation subject to the boundary and initial conditions in (\ref{heat ex 1 eq}). Here we use $\phi(r)=\sqrt{3/ \pi}e^{-3r^2}$ as our RBF, and choose $\delta=0.2$ and use points of the form $\frac{\mathbf{j}}{l}=(\frac{k}{l},\frac{m}{l})$ in $\Omega_{\delta}$, $k,m \in \mathbb{Z}$, $l=14$, to get the approximate solutions from Equation (\ref{nonhomohelm approx R2}). Next, the methods in \textsection{\ref{helm sec}} are used to solve problem (\ref{homo-helm a})-(\ref{homo-helm b}). $40$ collocation and source points are chosen in the same manner as in Example \ref{heatex} (see Figure \ref{heatfdm fig}). Finally, the solution, $\widetilde{u}_n(x,y)$ to the  problem (\ref{fdm helm eq}) subject to the initial and boundary conditions in (\ref{heat ex 1 eq}) is given by Equation (\ref{non homo helm approx sol}). By repeating this process in an iterative manner we can find $\widetilde{u}_2(x,y)$, $\widetilde{u}_3(x,y)$, $\cdots$, $\widetilde{u}_M(x,y)$.

  We evaluate the numerical solution, $\widetilde{u}_M$, and the true solution, $u$, at the same $6$ interior points $(x,y) \in \Omega$ and $T$ value used in Example \ref{heatex}. The results are presented below for different values of $M$.

 \begin{table}[H] 
 \centering
 \caption{True and numerical solutions to Example \ref{heatfdm} with $M=10$}
\begin{tabular}{ |p{3cm}||p{3cm}|p{3cm}|p{3cm}|  } 
 \hline
 $(x,y)\in \Omega$ & True \newline Solution & Numerical \newline Solution & Error  \\
 \hline
 $(-0.01,0.07)$ &   0.7014      & 0.7211 & 0.0197 \\ 
   $(-0.01,0.04)$ & 0.7682  &  0.7666& 1.5999e-03 \\
 $(-0.01,0.01)$ &  0.8008  &  0.7976  & 3.1984e-03  \\ 
  $(-0.01,-0.01)$ &   0.8008    & 0.8121 & 0.0113 \\
   $(-0.01,-0.04)$ &   0.7682   & 0.7721 &  3.9114e-03 \\
    $(-0.01,-0.07)$& 0.7014 &   0.7112 & 9.8114e-03\\
  
 \hline   
\end{tabular} 
\end{table}

  \begin{table}[H] 
 \centering
 \caption{True and numerical solutions to Example \ref{heatfdm} with $M=30$}
\begin{tabular}{ |p{3cm}||p{3cm}|p{3cm}|p{3cm}|  } 
 \hline
 $(x,y)\in \Omega$ & True \newline Solution & Numerical \newline Solution & Error  \\
 \hline
 $(-0.01,0.07)$ &   0.7014  & 0.7112 & 9.8321e-03\\ 
   $(-0.01,0.04)$ & 0.7682  & 0.7567 & 0.0114\\
 $(-0.01,0.01)$ &  0.8008  &  0.8100  &  9.2124e-03 \\ 
  $(-0.01,-0.01)$ &   0.8008    & 0.8025 &  1.7231e-03\\
   $(-0.01,-0.04)$ &   0.7682   & 0.7694  &  1.2234e-03 \\
    $(-0.01,-0.07)$& 0.7014 &     0.7100 & 8.6241e-03 \\
  
 \hline   
\end{tabular} 
\end{table}

    \begin{figure}[H]
\centering
  \includegraphics[width=.6\textwidth]{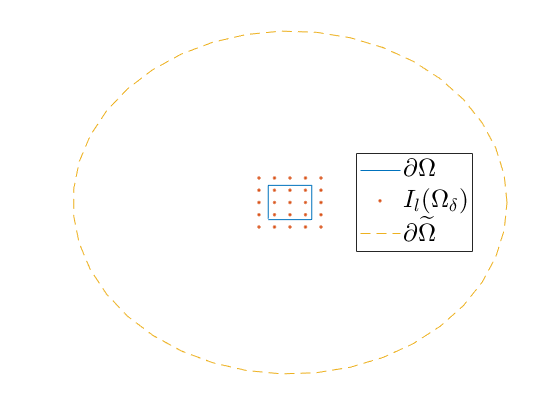} 
  \caption{Domains and points in $I_l(\Omega_{\delta})$ for Example \ref{heatfdm}} \label{heatfdm fig}
\end{figure}
  
  \end{exmp}

Next, we look at a wave problem using the Laplace transform.

\begin{exmp} \label{wave ex1}
Consider the following model of a vibrating rectangular membrane.
\begin{equation}
\begin{aligned}
    \frac{\partial^2 u}{\partial t^2}(x,y,t) &= \Delta u(x,y,t), && (x,y) \in \Omega, \quad & t > 0, \\
    u(x,y,t) &= 0, && (x,y) \in \partial \Omega, \quad & t > 0, \\
    \frac{\partial u}{\partial t}(x,y,0) &= xy, && (x,y) \in \Omega, \quad & t = 0, \\
    u(x,y,0) &= 0, && (x,y) \in \Omega, \quad & t = 0.
\end{aligned} \label{wave ex1 eq}
\end{equation}

  where $\Omega=(0,1) \times (0,1)$. Using the separation of variables method we find that the analytical solution to this problem is given by
  \begin{equation}
      u_e(x,y,t)=\sum_{m=1}^{\infty} \sum_{n=1}^{\infty}  \frac{4(-1)^{m+n}}{\pi^2  mn \lambda_{m,n}} \sin{(\lambda_{m,n}t)} \sin{(m \pi x)} \sin{(n \pi y)}, \label{wave ex1 analytic}
  \end{equation}
  where 
  $$
  \lambda_{m,n}=\pi \sqrt{m^2+n^2}.
  $$
  After applying the Laplace transform to the above problem, we choose $10$ fixed values $s_l$, $1 \leq l \leq 10$, to be the Laplace transform parameters mentioned in Equation (\ref{helm s value}), and then we obtain the following Helmholtz BVP
   \begin{alignat}{2}
     \Delta U(\mathbf{x})-s_l^2 U(\mathbf{x}) &= xy, \quad &&\mathbf{x} \in \Omega,  \label{wave ex 1 LT eq 1} \\
    U(\mathbf{x})&=0, \quad  && \mathbf{x} \in \partial \Omega, \label{wave ex 1 LT eq 2}
  \end{alignat}
  where $U(\mathbf{x})=\mathscr{L}[u(x,y,t)](s_l)$. For each $s_l$, problem (\ref{wave ex 1 LT eq 1})-(\ref{wave ex 1 LT eq 2}) can be solved by the methods in \textsection{\ref{helm sec}}. We will use $\phi(r)=\sqrt{\frac{10}{\pi}}e^{-10r^2}$ as our RBF, $\delta=0.1$, and $n=16$. Then we use points of the form $\frac{\mathbf{j}}{n}=(\frac{k}{n},\frac{m}{n})$ in $\Omega_{\delta}$, $k,m \in \mathbb{Z}$ to get the approximate solutions from Equation (\ref{nonhomohelm approx R2}). Next, we choose $N$ collocation points to be equally spaced points on $\partial \Omega$ and $N$ source points to be equally spaced points on a circle of radius $1.2$ whose center is $(0.5,0.5)$. Then the numerical solution to the problem (\ref{wave ex 1 LT eq 1})-(\ref{wave ex 1 LT eq 2}) is given by Equation (\ref{non homo helm approx sol}). Lastly, to find the numerical solution to the original problem (\ref{wave ex1 eq}) we use the formula for the inverse Laplace transform as mentioned above.
  \begin{figure}[H]
\centering
  \includegraphics[width=.6\textwidth]{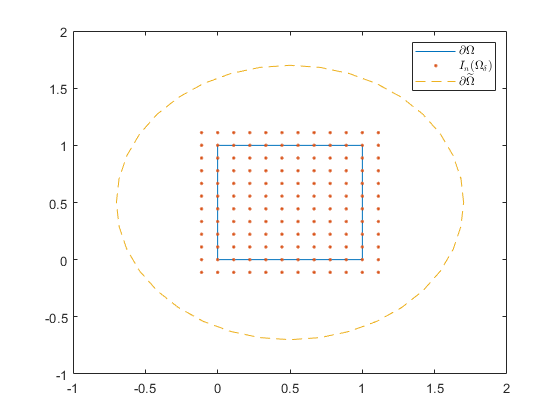}
  \caption{Domains and points used for Example \ref{wave ex1}}
\end{figure}
  
  To estimate the numerical error we will choose $6$ points $(x,y) \in \Omega$, and a time $t=T$, then calculate
  $$
  \text{error}=\vert u(x,y,T)-u_e(x,y,T) \vert.
  $$
  The errors are presented below for different values of $N$ and $T$.

         \begin{table}[H] 
 \centering
 \caption{Numerical errors to Example \ref{wave ex1} with $N=29$}
\begin{tabular}{ |p{3cm}||p{3cm}|p{3cm}|  } 
 \hline
 $(x,y)\in \Omega$ & T= 15& T=20   \\
 \hline
 $(0.14,0.35)$ & 7.8634e-03 & 2.0180e-03 \\ 
   $(0.29,0.35)$ & 6.0318e-03 &  0.0266  \\
 $(0.43,0.35)$ & 7.7390e-03 & 0.0647    \\ 
  $(0.57,0.35)$ &  0.0424    &  0.0523    \\
   $(0.71, 0.35)$ & 0.0728  &  0.0343      \\
    $(0.86, 0.35)$ & 0.0586  &  0.0212      \\
  
 \hline   
\end{tabular} 
\end{table}

       \begin{table}[H] 
 \centering
 \caption{Numerical errors to Example \ref{wave ex1} with $N=37$}
\begin{tabular}{ |p{3cm}||p{3cm}|p{3cm}|  } 
 \hline
 $(x,y)\in \Omega$ & T= 15& T=20   \\
 \hline
 $(0.14,0.35)$ & 9.0618e-04  &4.2558e-03  \\ 
   $(0.29,0.35)$ & 0.0257  & 0.0398    \\
 $(0.43,0.35)$ & 0.0180  & 0.0819   \\ 
  $(0.57,0.35)$ &  0.0446    &   0.0828  \\
   $(0.71, 0.35)$ &  0.0716  & 0.0732      \\
    $(0.86, 0.35)$ & 0.0547   & 0.0416  \\
  
 \hline   
\end{tabular} 
\end{table}

  \end{exmp}

\begin{exmp} \label{wave ex2}
  Let us revisit problem (\ref{wave ex1 eq}), but here we use the difference in time approach. Recall the analytical solution to this problem is given by Equation (\ref{wave ex1 analytic}).
  
  First we choose $M$ equally spaced points on the time interval $[0,T]$ for some $T$. Then we use the methods in \textsection{\ref{helm sec}} to solve Equation (\ref{time diff wave eq}) subject to the initial conditions (\ref{time diff wave ic}) and boundary condition (\ref{time diff wave bc}).  We use $\phi(r)=\sqrt{3/ \pi}e^{-3r^2}$ as our RBF.  We choose $\delta=0.2$ and use points of the form $\frac{\mathbf{j}}{l}=(\frac{k}{l},\frac{m}{l})$ in $\Omega_{\delta}$, $k,m \in \mathbb{Z}$, $l=14$, to get the approximate solutions from Equation (\ref{nonhomohelm approx R2}). Next, the methods in \textsection{\ref{helm sec}} are used to solve problem (\ref{homo-helm a})-(\ref{homo-helm b}). $N$ equally spaced points on $\partial \Omega$ are used as the collocation points. The source points are $N$ equally spaced points on a circle of radius $1.2$ with center $(0.5,0.5)$. Finally, the solution, ${u}_n(x,y)$ to the  problem (\ref{time diff wave eq}) with initial conditions (\ref{time diff wave ic}) and boundary condition (\ref{time diff wave bc}) is given by Equation (\ref{non homo helm approx sol}). Then repeating this process in an iterative manner we can find $u_3(x,y)$, $u_4(x,y)$, $\cdots$, $u_M(x,y)$.
  
  To get the numerical error, we choose $6$ points $(x,y)\in \Omega$ and calculate
  $$
  \text{error}= |u_e(x,y,T)-u_M(x,y,T)|,
  $$
  where $u_e$ is the exact solution. The results are presented below for different values of $T$, $M$, and $N$.

       \begin{table}[H] 
 \centering
 \caption{Numerical errors to Example \ref{wave ex2} with $N=21$ and $M=15$}
\begin{tabular}{ |p{3cm}||p{3cm}|p{3cm}|  } 
 \hline
 $(x,y)\in \Omega$ & T= 5& T=15   \\
 \hline
$(0.14,0.35)$  & 0.0197   &  0.0138  \\ 
$(0.29,0.35)$  & 0.0151  &  5.0604e-05  \\
$(0.43,0.35)$  &  6.1836e-03  &   1.0821e-03  \\ 
$(0.57,0.35)$  &  0.0377 &  0.0343  \\
$(0.71, 0.35)$ &  0.0260   &  0.0615  \\
$(0.86, 0.35)$ &  8.3917e-03  &  0.0350  \\
  
 \hline   
\end{tabular} 
\end{table}
   
       \begin{table}[H] 
 \centering
 \caption{Numerical errors to Example \ref{wave ex2} with $N=157$ and $M=10$}
\begin{tabular}{ |p{3cm}||p{3cm}|p{3cm}|  } 
 \hline
 $(x,y)\in \Omega$ & T= 5& T=15   \\
 \hline
$(0.14,0.35)$  &  7.4234e-04  & 0.0146   \\ 
$(0.29,0.35)$  &  0.0239   &  8.7654e-04  \\
$(0.43,0.35)$  &  0.0171  &   8.0491e-04  \\ 
$(0.57,0.35)$  &  3.0332e-03  &  0.0355   \\
$(0.71, 0.35)$ &  9.9534e-03  &  0.0660   \\
$(0.86, 0.35)$ &   7.9493e-03 &  0.0521   \\
  
 \hline   
\end{tabular} 
\end{table}

     \begin{table}[H] 
 \centering
 \caption{Numerical errors to Example \ref{wave ex2} with $N=157$ and $M=15$}
\begin{tabular}{ |p{3cm}||p{3cm}|p{3cm}|  } 
 \hline
 $(x,y)\in \Omega$ & T= 5& T=15   \\
 \hline
$(0.14,0.35)$  &1.1104e-04 &  8.8195e-04  \\ 
$(0.29,0.35)$  &  0.0225  &   0.0214  \\
$(0.43,0.35)$  &  0.0155  &  0.0143   \\ 
$(0.57,0.35)$  &  1.3922e-03 &  1.4931e-04   \\
$(0.71, 0.35)$ &  0.0113  &  0.0124  \\
$(0.86, 0.35)$ &  8.8048e-03  &  9.5762e-03  \\
  
 \hline   
\end{tabular} 
\end{table}

  \end{exmp}

\section{Conclusions}
We investigated diffusion and wave initial boundary value problems by transforming the respective problems into Helmholtz boundary value problems and then using a meshless approach to solve the problems. One of the motivations for using meshless methods is its simplicity in implementation. As shown in the examples above, many of the problems can be solved using relatively few points.   

\bibliographystyle{plain}
%\bibliography{citations.bib}

\end{document}